\documentclass[12pt]{article}
\newtheorem{prethm}{{\bf Theorem} {\rm (MAIN)}}

\newtheorem{precor}{{\bf Corollary}}

\newenvironment{corollary}{\begin{precor}{\hspace{-0.5
               em}{\bf.}}}{\end{precor}}
\newtheorem{preprop}{{\bf Proposition}}

\newtheorem{preprob}{{\bf Problem}}

\newenvironment{problem}{\begin{preprob}{\hspace{-0.5
               em}{\bf.}}}{\end{preprob}}
\newtheorem{preque}{{\bf Question}}

\newenvironment{question}{\begin{preque}{\hspace{-0.5
               em}{\bf.}}}{\end{preque}}
\newtheorem{prelemma}{{\bf Lemma}}

\newenvironment{lemma}{\begin{prelemma}{\hspace{-0.5
               em}{\bf.}}}{\end{prelemma}}
\newtheorem{preex}{{\bf Example}}

\newenvironment{example}{\begin{preex}{\hspace{-0.5
               em}{\bf.}}}{\end{preex}}
\newtheorem{prepro}{{\bf Proposition}}

\newenvironment{proof}{{\bf Proof.}}{\hfill\rule{2mm}{2mm}}
\newtheorem{predefinition}{Definition}

\newenvironment{definition}{\begin{predefinition}{\hspace{-.7
               em}{\bf.}}}{\end{predefinition}}
\newtheorem{pretheorem}{{\bf Theorem}}

\newenvironment{theorem}{\begin{pretheorem}{\hspace{-0.5
               em}{\bf.}}}{\end{pretheorem}}
\def\newpic#1{}
\def\NP{{\bf NP}}
\def\spec{{\rm Spec}}
\def\bracket{{\rm bracket}}
\def\sbracket{{\rm skew \ bracket}}
\def\2switch{{\rm matching 2--switch}}
\def\comment#1{}
\newtheorem{prelem}{{\bf Theorem}}

\newenvironment{lem}{\begin{prelem}{\hspace{-0.5
               em}{\bf.}}}{\end{prelem}}
\title{\bf  On the spectrum of the forced matching number of  graphs
}

\author{
{\bf P. Afshani$^a$, H. Hatami$^a$, and E.S.
Mahmoodian$^b$} \\
 \\[1mm]
{\small Institute for Studies in Theoretical Physics and Mathematics (IPM) }\\
 {and}
\\$^a${\small\it Department of Computer Engineering} \\
$^b${\small\it Department of Mathematical Sciences}\\
Sharif University of Technology \\
P.O. Box 11365--9415, Tehran, I.R. Iran}
\begin{document}
\maketitle

\begin{abstract}
Let $G$ be a graph that admits a perfect matching. A {\sf forcing
set} for a perfect matching $M$ of $G$ is a subset $S$ of $M$,
such that $S$ is contained in no other perfect matching of $G$.
This notion originally arose in chemistry in the study of
molecular resonance structures. Similar concepts have been studied
for block designs and graph colorings under the name {\sf defining
set}, and for Latin squares under the name {\sf critical set}.
Recently several papers have appeared on the study of forcing sets
for other graph theoretic concepts such as dominating sets,
orientations, and geodetics. Whilst there has been some study of
forcing sets of matchings of hexagonal systems in the context of
chemistry, only a few other classes of graphs have been
considered.

Here we study the spectrum of possible forced matching numbers for
the grids $P_m \times P_n$, discuss the concept of a forcing set
for some other specific classes of graphs, and  show that the
problem of finding the smallest forcing number of graphs is
\NP--complete.
\end{abstract}
\begin{center}
\begin{figure}[ht]
\input{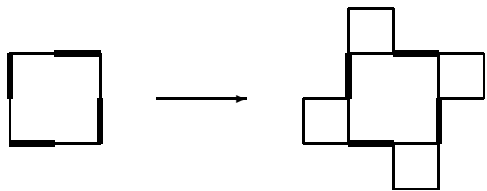}
\vspace*{-1cm} \caption{The graph $G_4$, obtained from $C_8$.
\label{G_4}}

\end{figure}
\end{center}

\newline
{{\sc Keywords:} forcing number; matching in graphs; spectrum;
hypercubes; grids.

\section{Introduction  and preliminaries}
Let $G$ be a graph that admits a perfect matching. A {\sf forcing
set} for a perfect matching $M$ of $G$ is a subset $S$ of $M$,
such that $S$ is contained in no other perfect matching of $G$.
\begin{example}
In Figure~{\rm\ref{m_s}} a forcing set of size $6$ is shown for a
matching in an $8 \times 12$ grid that is $P_8 \times P_{12}$. The
bold edges form a matching, and the edges in the forcing set are
indicated by small circles.
\begin{center}
\begin{figure}[ht]
\input{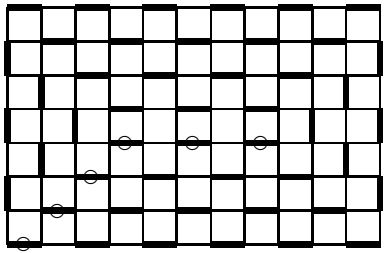}
\vspace*{1mm} \caption{\label{m_s}A forcing set for a matching in
$P_8 \times P_{12}$.}
\end{figure}
\end{center}
\end{example}
The matching in the Example~\ref{m_s} has a pattern which will be
used in the next section. It is called a {\sf concentrated
alternating cycles matching} or a  {\sf CACM} of size $8\times
12$, and is defined in general for a $P_{2m}\times P_{2n}$ as
follows: a CACM of size $2m \times 2n$ is a special matching in
$P_{2m}\times P_{2n}$, in which the vertices of the first row and
also the last row are matched horizontally, and the remaining
vertices of the first column and the last column are matched
vertically, so that these matched edges form an alternating cycle.
We continue this process recursively for  the remaining vertices,
which form a grid of size $(2m-2)\times (2n-2)$.

The cardinality of a smallest forcing set of $M$ is called {\sf
forced matching number}, and is denoted by $f(G,M)$, which we will
henceforth call the {\sf forcing number of $M$}. Also $f(G)$ and
$F(G)$, respectively, denote the minimum and maximum of $f(G,M)$
over the set of all perfect matchings $M$ of $G$. As all our
matchings will be perfect, we drop the use of ``perfect'' after
this point.

The notion of a forcing number originally arose in chemistry in
1987 in the study of molecular resonance
structures~\cite{MR88d:05159}. Later, in~\cite{MR94b:05193},
Harary introduced the concept of the forcing number of a perfect
matching and of other concepts in graphs. Since then, papers have
appeared on the forced orientation number of
 graphs~\cite{MR96m:05094,FarzadMahdianMahmoodianetal}
, dominating sets~\cite{MR98f:05083}, and
geodetics~\cite{MR2000d:05037}.

 Similar concepts have been studied under the
name {\sf defining set} for block
designs~\cite{MR91e:05016,MR96h:05021} for graph
colorings~\cite{MR98b:05044}, and under the name {\sf critical
set} for Latin squares~\cite{MR80j:05022,MR2000g:05034}. There
has been some study of forcing sets of matchings of hexagonal
systems (in the context of chemistry), and only a few other
classes of graphs have been
considered~\cite{MR1105504,PachtersThesis,MR99c:05156,StopSign,Riddle}.
One of the interesting problems is the study of the spectrum of
forcing numbers of a given graph; to this end, the following
definition is taken from~\cite{AdamsMahdianMahmoodian}.
\begin{definition}
The {\sf spectrum of forcing numbers} for a graph $G$ is a set of
natural numbers defined as: \\
$\spec (G)=\{ k \ | \ $ there exists a matching $M$ of $G$ such
that $f(G,M)=k\}$.
\end{definition}
The spectra of hypercubes is studied
in~\cite{AdamsMahdianMahmoodian}. In Section~2, we study the
spectrum of $P_m \times P_n$ and show that there are no gaps in
the spectra of forcing numbers of certain types of graphs which
include $P_m \times P_n$ and stop signs. In Section~3, we further
discuss the concept of forcing numbers for some specific classes
of graphs such as $P_m \times P_n$, $C_m \times P_n$, and $C_{2n}
\times C_{2n}$. Finally in Section~4, we investigate the
computational complexity of the problem of finding the forcing
number of a graph.
%
\section{Spectrum}
A natural question is: Which finite subsets of natural numbers are
the spectra of some graph or other? In order to answer this
question we need the following lemma.
\begin{lemma}
\label{spcH}
  If $G$ is a graph with $\spec(G)=A$, then for any integer $k$, there exists
  a graph $H$ with $\spec (H)=\{x+k \ | \ x \in A\}$.
\end{lemma}
\begin{proof}
The graph $H$ can be constructed by adding a union of $k$ disjoint
copies of $C_4$ (cycles of size 4) to $G$. Trivially $\spec
(H)=\{x+k \ | \ x \in A\}$.
\end{proof}

Next, for a given $n$ we define a graph $G_n$ by replacing every
other edge in $C_{2n}$ by a cycle of size $4$. This is illustrated
for $n=4$ in Figure~\ref{G_4}. Any of the bold edges from $C_8$
forces a matching in $G_4$. These edges are called {\sf forcing
edges}}.
The following trivial lemma is to facilitate the proof of the
subsequent theorem .
\begin{lemma}
\label{SpecGi} We have: \ \
 $\spec (G_n)=\{1,n\}$.
\end{lemma}
\begin{theorem}
  For any finite set $A \subset N$, there exists a
  graph $G$ with $\spec (G)=A$. Indeed, $G$ can be chosen to be a planar
  bipartite graph.
\end{theorem}
\begin{proof} Using Lemma~\ref{spcH}, we can assume that $1\in A$.
Firstly, for every $i \in A$ ($i \neq 1$), we assign a
corresponding graph $G_i$, each of which has some edges which are
forcing edges. Construct a graph $G$ by ``gluing'' each of these
$G_i$ to a common {\it forcing edge} $e=\{u,v\}$. We claim that
$\spec (G)=A$. Indeed $\{ e \}$ is a forcing set for $G$. Thus $1
\in \spec (G)$. Now, if we have a matching $M$ which {\it does
not} contain $e$, then both ends of $e$ must be matched with some
other vertices in one of the $G_i$, say $G_l$. Then $M$ generates
$l$ disjoint alternating cycles of size $4$ in $G_l$, so any
forcing set of $M$ has at least $l$ edges from $G_l$. Also observe
that a forcing set of size $l$ for $G_l$ is also a forcing set for
$G$. In fact the constructed graph $G$ is planar and bipartite.
\end{proof}

Next, we study the spectra of some special graphs. First we give a
simple proof of a theorem determining the spectrum of the grid
$P_{2n} \times P_{2n}$. We then generalize that proof, to show
that there are no gaps in the spectra of some specific graphs
including $P_{2m} \times P_{2n}$ and stop signs. Recall that an
{\sf $(n,k)$ stop sign} ($k \le n-1$) is a graph obtained from
$P_{2n} \times P_{2n}$ by deleting all of the vertices along the
$k$ diagonals closest to each of the four corners~\cite{StopSign}.

 So our result is that the spectrum of any
such graph contains all the numbers between the smallest and the
largest forcing number. Hence if we find the largest and the
smallest forcing number for those graphs, then the spectrum is
precisely determined.
%
\begin{definition}
    A {\sf matching 2--switch} is an operation on a graph defined
    by the replacement of matching edges with nonmatching edges in
     an alternating cycle of
    size four.
\end{definition}
\noindent The following lemma and its immediate corollary are
instrumental to our results.
\begin{lemma}
    A \2switch on a matching $M$ does not change the forcing
    number by more than $1$.
\end{lemma}
\begin{proof}
Suppose that $e_1=\{u_1,v_1\}$ and $e_2=\{u_2,v_2\}$ are two edges
of $M$ that form an alternating cycle ($u_1 v_1 v_2 u_2$). At
least one of these two edges must be in the forcing set of $M$.
Now consider a new matching $M'$ which is obtained by removing the
edges $e_1$ and $e_2$ from $M$, and adding $e'_1=\{u_1,u_2\}$ and
$e'_2=\{v_1,v_2\}$ to it. If $S$ is a forcing set for $M$, then
$(S \cup \{ e'_1 , e'_2\}) \backslash \{ e_1 , e_2 \}$ is a
forcing set for $M'$, so the forcing number of $M'$ is at most one
more than the forcing number of $M$. The same argument holds when
we convert $M'$ to $M$.
\end{proof}

\begin {corollary}
\label{2switch}
       In a graph $G$
       with a sequence of matchings
       $M_1, M_2, \ldots, M_s$, such that
       $M_{i+1}$ is obtained from $M_i$ by a \2switch,
       all the numbers between $f(G,M_1)$ and $f(G,M_s)$ appear in the
       set consisting of the
       forcing numbers of $M_1, M_2, \ldots, M_s$.
\end{corollary}
Now we are ready to determine the spectrum of forcing numbers of
$P_{2n} \times P_{2n}$. Pachter and Kim proved the following
theorem.
\begin{lem}
\label{patcherp2np2n}{\rm\cite{MR99c:05156}}
 Let $M$ be a
matching of $P_{2n} \times P_{2n}$. Then\ \ $n \leq f(P_{2n}
\times P_{2n},M) \leq n^2.$
\end{lem}
In the following theorem we show that $f(P_{2n}\times P_{2n},M)$
actually takes on {\it all} the values between $n$ and $n^2$.
\begin{theorem}
\label{specp2np2n} We have: \ \ $\spec (P_{2n} \times P_{2n})=\{n
, \ldots, n^2\}$.
\end{theorem}
\begin{proof}
By Corollary~\ref{2switch}, it is sufficient to convert a matching
with forcing number $n^2$ to a matching with forcing number $n$,
by repeatedly applying matching 2--switches. We illustrate a
process for this, using the example graph $P_6 \times P_6$ (that
is when $n=3$) in Figure~\ref{process}. The matching $M_1$ has
forcing number $n^2$ (=9), and $M_s$ which is a CACM has forcing
number $n$ (=3). It is easily seen that it is possible to convert
$M_1$ to $M_2$ and $M_2$ to $M_3$ by applying matching
2--switches. By performing the same operations recursively on the
inner $(2n-2)\times (2n-2)$ grid in $M_3$, we finally obtain the
matching $M_s$.
\begin{center}
\begin{figure}[ht]

\input{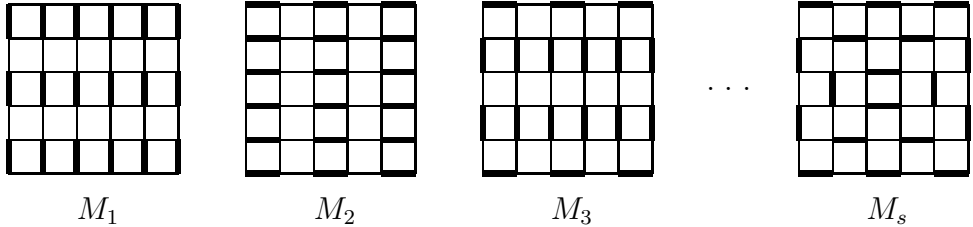}
\vspace*{-5mm}
 \caption{Applying matching $2$--switches to reduce forcing numbers.\label{process}}
\end{figure}
\end{center}
It should be easy to see that this procedure is valid for any $n$.
Since $f(P_{2n}\times P_{2n},M_1) = n^2$ and $f(P_{2n}\times
P_{2n},M_s) = n$, so $\spec (P_{2n} \times P_{2n})=\{n , \ldots,
n^2\}$.
\end{proof}

\begin{center}
\begin{figure}[ht]
\input{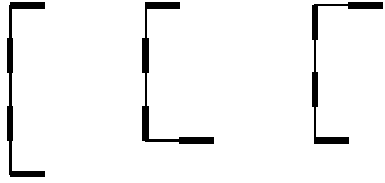}
\vspace*{-1cm} \caption{An $(i,j,2)- \bracket$ and $(i,j,2)$-skew
brackets of type I and II.\label{bracket}}
\end{figure}
\end{center}

Next we generalize the method applied in the proof of
Theorem~\ref{specp2np2n} to more general graphs. To facilitate
this, we label the vertices of $P_n \times P_n$ by ordered pairs
$(i,j)$, where $ 1 \leq i, j \leq n$; and $i$ is the row number
and $j$ is the column number of that vertex.
\begin{definition}
An induced subgraph $G$ of a grid with vertex set $V(G)$ is called
a {\sf column continuous subgrid} if it has the following
property:
\begin{itemize}
\item If $(i_1,j), (i_2,j) \in V(G)$, then for all integers
 $i$, such that $i_1 \leq i \leq i_2$,
we have $(i,j) \in V(G)$.
\end{itemize}
\end{definition}
Suppose $G$ is an induced subgraph of $P_n \times P_n$ which has a
matching $M$. An $(i,j,k)-${\sf bracket} is a bracket shaped
subset of the edges of $M$ (e.g. Figure~\ref{bracket}) as in the
following:

$
\begin{array}{l}
 \{ \ \ \{(i,j),(i,j+1)\},\\
 \ \ \ \{(i+1,j),(i+2,j)\},\{(i+3,j),(i+4,j)\},
\ldots,  \{(i+2k-1,j),(i+2k,j)\}, \\
\ \ \ \{(i+2k+1,j),(i+2k+1,j+1)\}\ \ \}; \\
\end{array}
$

 and the following set of edges is called an $(i,j,k)-${\sf skew bracket (of type I)} $(k > 0)$:

$
\begin{array}{l}
 \{  \ \ \{(i,j),(i,j+1)\},\\
  \ \ \ \{(i+1,j),(i+2,j)\},\{(i+3,j),(i+4,j)\}, \ldots,
 \{(i+2k-1,j),(i+2k,j)\}, \\
 \ \ \ \{(i+2k,j+1),(i+2k,j+2)\}\ \ \}.\\
\end{array}
$

{\sf Skew bracket (of type II)} is defined similarly as the
following set of edges:

$
\begin{array}{l}
 \{  \ \ \{(i,j+1),(i,j+2)\},\\
  \ \ \ \{(i,j),(i+1,j)\},\{(i+2,j),(i+3,j)\}, \ldots,
 \{(i+2k-2,j),(i+2k-1,j)\}, \\
 \ \ \ \{(i+2k,j),(i+2k,j+1)\}\ \ \}.\\
\end{array}
$

See Figure~\ref{bracket} for an example.
\begin{lemma}
\label{brack1} Let $G$ be a column continuous subgrid of $P_n
\times P_n$. If $M$ is a matching in $G$ which contains an
$(i,j,k)-\bracket$, then we can apply matching $2$--switches to
$M$ on all the edges which have both endpoints in the following
set of vertices:
\[ \{ (a,b) \ | \ i \leq a \leq i+2k+1, j \leq b \leq n \} \cap V(G), \]
so that the resulting matching contains the following edges:
\[ \{\{(i,j),(i+1,j)\},\{(i+2,j),(i+3,j)\}, \ldots,
\{(i+2k,j),(i+2k+1,j)\}\}.
\]
\end{lemma}
\begin{proof}
Note that we want to show that $M$ can be changed to a matching
such that all the edges in it which touch the set of vertices
$(a,j)$ in the $j$-th column, for $i \leq a \leq i+2k+1$, are all
vertical. We apply mathematical induction on $k$. The case $k=0$
is trivial. Suppose the statement is true for $k=p$. Consider an
$(i,j,p+1)-\bracket$. There are two cases.

The first case is where all the edges of $M$ which touch the set
of vertices $A=\{(i+1,j+1), (i+2,j+1), \ldots, (i+2p+2,j+1)\}$ are
all vertical (obviously $ A \subseteq V(G)$). It is easy to verify
the lemma in this case.

If it is not the first case, then some of the edges which touch
the set $A$ are horizontal. The horizontal and vertical edges
which touch $A$ make some $(x,j+1,t)-\bracket$s. We choose one of
these brackets and apply the induction hypothesis to it,
increasing the number of vertical edges which touch $A$ by $1$. By
repeating this process we can convert all of the matching edges
touching $A$ to vertical matching edges, which is the first case.
Note that the induction hypothesis ensures that converting an
$(x,j+1,t)-\bracket$ does not have any effect on previously
converted vertical edges.
\end{proof}

\begin{corollary}
\label{nobrack} Let $G$ be a column continuous subgrid. If $M$ is
a matching in $G$, then by applying matching $2$--switches we can
convert $M$ to a matching which contains no $(i,j,k)-\bracket$.
\end{corollary}
\begin{proof}
Let $j$ $(1 \leq j \leq n)$ be the minimum value for which there
exists some bracket in the $j$-th column. By using
Lemma~\ref{brack1}, we can destroy this bracket by matching
2--switches. If we continue this process, there will be no
bracket left in this column, and so the value of $j$ increases.
Repeating this process removes all brackets.
\end{proof}

\begin{lemma}
\label{noskew}
 Let $G$ be a column continuous subgrid. If $M$ is a
matching in $G$ in which there is no bracket, then there is also
no skew bracket of any type in $M$.
\end{lemma}
\begin{proof}
Assume to the contrary that $M$ has no bracket, but that there
does exist for example an $(i,j,k)-\sbracket$ of type I in $M$.
Since there are odd number of vertices in the set $\{(i+1,j+1),
(i+2,j+1), \ldots,(i+2k-1,j+1)\}$, the presence of matching edges
in the $(i,j,k)-\sbracket$ leads to the presence of at least one
bracket in the column $j+1$. A contradiction. similar argument
holds, if we assume that $M$ contains a skew bracket of type II.
\end{proof}

\begin{theorem}
\label{nogap} There are no gaps in the spectrum of a column
continuous subgrid.
\end{theorem}
\begin{proof}
Assume that $G$ is a column continuous subgrid. We show that it is
possible to convert a given matching of $G$ to any other matching
of $G$, by applying matching 2--switches.

Suppose we have two matchings in $G$. By Corollary~\ref{nobrack}
we remove all brackets from both of these matchings and end up
with matchings say $M$ and $M'$. If $M \neq M'$, then there exists
a cycle which is alternating in $M$ and $M'$. So if we consider
the first column which is touched by this cycle, at least one of
$M$ and $M'$ contains either a bracket or a skew bracket, and this
contradicts Lemma~\ref{noskew} for neither $M$ nor $M'$ contains a
bracket.
\end{proof}

Note that the assumption that the graph involved is an ``induced
subgraph'' of a grid is necessary for the result of
Theorem~\ref{nogap}. Also the assumption that it be ``column
continuous'' is necessary, as can be seen from the fact that
$\spec (G_4)=\{1, 4\}$, where $G_4$ is shown in Figure~\ref{G_4}.
Indeed one can give infinitely many examples to show the
necessity of this condition.

Since both $P_m \times P_n$ and the $(n,k)$ stop sign are column
continuous subgrids, we have the following corollary.
\begin{corollary}
\label{gappmpn}
 There are no gaps in the spectrum of forcing numbers of $P_m \times P_n$
 and in the spectrum of forcing numbers of an $(n,k)$ stop
sign.
\end{corollary}
The spectra of stop signs follow from the following theorem and
Corollary~\ref{gappmpn}.
\begin{lem}
\label{stop}{\rm\cite{StopSign}} Let $G$ be an $(n,k)$ stop sign
and $M$ be a matching of $G$. The forcing number of $M$ is
bounded by
\[ n \leq f(G,M) \leq (n-\lceil \frac{k-1}{2} \rceil)
(n-\lfloor \frac{k+1}{2} \rfloor),\] and the bounds are sharp.
\end{lem}
\section{Some special classes of graphs}
In this section we study $F(G)$, where $G$ is from some special
classes of graphs: a product of two paths, a product of a cycle
and a path, or a product of two  cycles. We also introduce an
upper bound for the smallest forcing number of a product of two
paths. Pachter and Kim pointed out the following useful result.
\begin{lem}
{\rm \cite{MR99c:05156,lucchesi}} \label{planar}
 If $G$ is a planar bipartite graph and $M$ is a
 matching in $G$, then the forcing number of $M$ is equal to
the maximum number of disjoint $M$-alternating cycles.
\end{lem}
\subsection{$P_{m}\times P_{n}$}

Applying the same method as in~\cite{MR99c:05156} we see that:
%
$$F(P_{m} \times P_{n}) = \lfloor \frac{m}{2} \rfloor \cdot \lfloor \frac{n}{2} \rfloor.$$
In contrast, finding $f(P_{m} \times P_{n})$ does not seem to be
so easy. We introduce a pattern which gives an upper bound for it.
\begin{theorem}
\label{pmpn} We have:
\begin{enumerate}
\item[{\bf{(i)}}]
$f( P_{2k} \times P_{(2k+1)l+r}) \leq kl+\lceil \frac{r-1}{2}
\rceil$,
 \ where $0 \leq r \leq 2k {\it \ and} \ l \ge 1;$
\item[{\bf{(ii)}}]
$f( P_{2k+1} \times P_{(2k+2)l+2r}) \leq kl+r$, \ where $0 \leq
2r \leq 2k+1 {\it \ and} \ l \ge 1.$
\end{enumerate}
\end{theorem}
\begin{proof}
We construct a matching $M$ for which there is a forcing set of
the desired size in the statement of the theorem.

\noindent {\bf  (i)} \ We choose the following $l$ columns: $1,
(2k+1)+1, 2(2k+1)+1, \dots, (l-1)(2k+1)+1$; and also the last
column if $r$ is even. There are $2k$ vertices in each column, we
take a matching in each of the chosen columns. Ignoring the chosen
columns we have $l$ blocks of size $2k \times 2k$ (this is
strictly true for all but the last block, which is not in general
square, but is of height $2k$ and of width varying with $r$). We
substitute a CACM of appropriate size into each one of these
blocks (see Figure~\ref{block}).

This matching $M$ has a forcing set of size $k(l-1)+ \lceil
\frac{2k+r-1}{2} \rceil = kl+ \lceil \frac{r-1}{2} \rceil$ as
shown in Figure~\ref{block}.
\begin{center}
\begin{figure}[ht]
\input{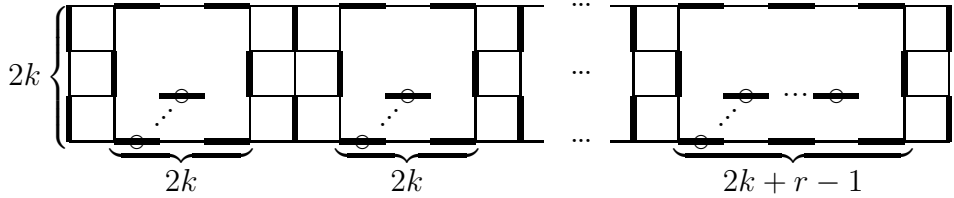}
%
 \caption{The pattern of $M$ when $r$ is odd.\label{block}}
\end{figure}
\end{center}
In the following figure $M$ is demonstrated for $P_8 \times
P_{25}$.
\begin{center}
\begin{figure}[ht]
\input{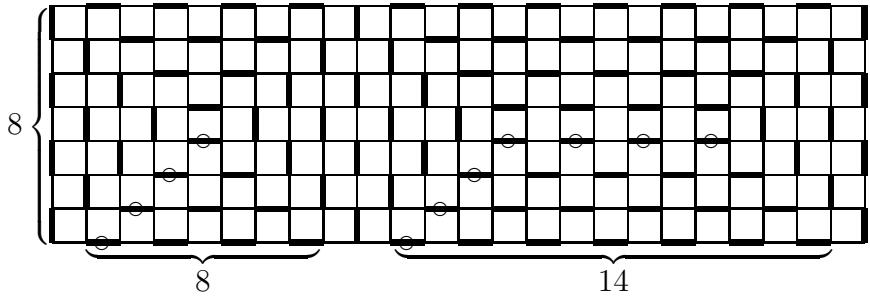}
\vspace*{3mm}\caption{A forcing set of size $11$ for $P_8 \times
P_{25}.$\label{8-25}}
\end{figure}
\end{center}
{\bf (ii)} To deal with this case we construct a matching in a
similar fashion to that of the previous case. To facilitate this,
we introduce some notation. A {\sf UCACM} and a {\sf DCACM} of
size $(2m-1) \times 2n$ are built from a CACM of size $2m \times
2n$ by removing the vertices of the first row, and the last row,
respectively.

In this case we partition $P_{2k+1} \times P_{(2k+2)l+2r}$ to
$(l-1)$ blocks of size $(2k+1)\times (2k+2)$ and one last block
of size $(2k+1)\times (2k+2r+2)$, and then replace each block
alternatively with a UCACM or a DCACM of appropriate size. This
is illustrated in Figure~\ref{p5-28} for the case $P_5 \times
P_{28}$.
\begin{center}
\begin{figure}[ht]
\input{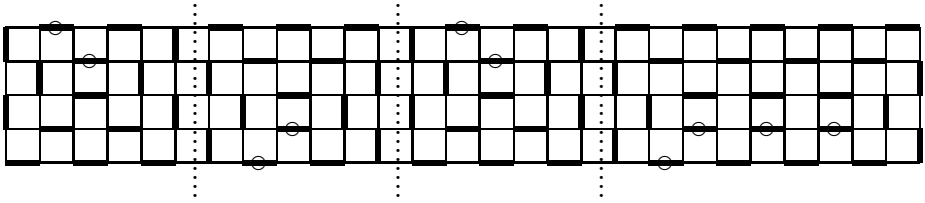}
\vspace*{1mm} \caption{A forcing set of size $10$ for $P_5 \times
P_{28}.$\label{p5-28}}
\end{figure}
\end{center}
Again the resulting matching has a forcing set of the desired
size.
\end{proof}

Note that in the previous theorem, in each case there are
appropriate number of $M$-alternating disjoint cycles which
Theorem~\ref{planar} implies that the size of the corresponding
forcing sets are smallest. Based on observations of small cases,
we conjecture that the bounds in Theorem~\ref{pmpn} are sharp.
%
\subsection{$P_{m} \times C_{2n}$}
The following theorem gives the exact value for the size of a
largest forcing set for $P_m \times C_{2n}$.
\begin{theorem}
\label{pmc2n} For every $k,n \geq 1$ we have: $$F(P_{2k} \times
C_{2n})=kn \ {\rm and} \  F(P_{2k+1} \times C_{2n})=kn+1.$$
\end{theorem}
\begin{proof}
Consider $P_{m} \times C_{2n}$ drawn as $2n$ ``vertical'' copies
of $P_{m}$ and $m$ ``horizontal'' copies of $C_{2n}$ on the set
of vertices in the columns. The graph $P_{m} \times C_{2n}$ is
planar and bipartite, so by Theorem~\ref{planar} for any matching
$M$, $f(P_{m} \times C_{2n}, M)$  is equal to the maximum number
of disjoint $M$-alternating cycles.

Since the girth of $P_{2k} \times C_{2n}$ is $4$, its largest
forcing number is not greater than $\frac{4kn}{4}=kn$. A matching
which has all edges horizontal clearly has forcing number $kn$.

For $P_{2k+1} \times C_{2n}$, suppose that $M$ is a matching, and
let ${\cal A}$ be a set of disjoint $M$-alternating cycles. If
there is an $M$-alternating cycle in ${\cal A}$ which intersects
a column exactly {\it once}, then it is at least of size $2n$. In
this case there are at most $\frac{(2k)(2n)}{4}=kn$ other cycles
in ${\cal A}$, and we are done.

So assume that there is no $M$-alternating cycle in ${\cal A}$
which intersects some column in exactly one vertex. In ${\cal A}$,
each cycle has at least two vertices of intersection with each
column that it intersects, so each column intersects at most $k$
cycles in ${\cal A}$. Now, as there are $2n$ columns if we count
all cycles, we get $k(2n)$. But in this way each cycle is counted
at least twice, as it intersects at least two {\it different}
columns. So there are at most $\frac{k(2n)}{2}=kn$ cycles.

In this case, a matching which has all edges horizontal clearly
has forcing number equal to $kn+1$.
\end{proof}

The following interesting problems remain open.
\begin{problem}
Find \ \ $F(P_{2m} \times C_{2n+1})$.
\end{problem}
\begin{problem}
Find \ \ $f(P_{m} \times C_{n})$.
\end{problem}
%
\subsection{$C_{2n} \times C_{2n}$}
It is conjectured in~\cite{Riddle} that $F(C_{2n} \times
C_{2n})=n^2$. A result in this direction is given in the following
theorem.
\begin{theorem}
We have: \ \ $F(C_{2n} \times C_{2n}) \leq n^2 +  \frac{n}{2}.$
\end{theorem}
\begin{proof}
Let $M$ be a  matching in  $C_{2n} \times C_{2n}$ which has the
largest forcing number. We show that there exists a forcing set of
size less than or equal to  $n^2 + \frac{n}{2}$ for $M$. The
number of edges in $M$ is $2n^2$, and at least $n^2$ of these
edges are in the same direction (``horizontal'' or ``vertical'').
Without loss of generality, suppose at least $n^2$ of the edges in
$M$ are horizontal. So there exists a row, say $r$ in which at
least $\frac{n}{2}$ edges of $M$ are horizontal. Thus, there are
at most $n+\frac{n}{2}$ matching edges which touch this row. We
take all these matching edges in our forcing set.

Removing the vertices we chose in our forcing set, we get a planar
graph, and we consider two cases. First, the case
 in which all the matching edges of row $r$ are horizontal.
 In this case, we have already chosen $n$ edges and the rest of the
 graph is a $P_{2n-1} \times C_{2n}$, which by Theorem~\ref{pmc2n}
 needs at most $n(n-1)+1$ edges to be forced. In the second case,
we have chosen at most $n+\frac{n}{2}$ edges and the graph
obtained after deleting those vertices has at most $2n-1$ vertices
in each column and also has at least one column with exactly
$2n-2$ vertices. Since we have a column which contains $2n-2$
vertices, by using the technique of the previous theorem, we can
say that the largest forcing number of the resulting graph is at
most $n(n-1)$.
So the forcing number of $M$ is at most
$n(n-1)+n+\frac{n}{2}=n^2+\frac{n}{2}$.
\end{proof}
%
\section{Computational complexity}
In~\cite{AdamsMahdianMahmoodian}, Adams, Mahdian, and Mahmoodian
studied the following problem and gave a proof for its
\NP-completeness.
\noindent
\begin{itemize}
\item
{\sc Smallest forcing set problem} \\ [1mm] {\sc Instance:} A
graph $G$, a matching $M$ in $G$, and an
integer $k$. \\
{\sc Question:} Is there any subset $S$ of at most $k$ edges of
$M$, such that $S$ is a forcing set for $M$?
\end{itemize}
\begin{lem}~{\rm\cite{AdamsMahdianMahmoodian}}
\label{AdamsMahdianMahmoodian}{\sc Smallest forcing set} is
\NP-complete for bipartite graphs with maximum degree $3$.
\end{lem}
They also left an open question which we answer in this section.
The question is finding the computational complexity of the
following problem:
\noindent
\begin{itemize}
\item
{\sc Smallest forcing number of graph} \\ [1mm] {\sc Instance:} A
graph $G$ and an
integer $k$. \\
{\sc Question:} Is there any matching in $G$ with the forcing
number of at most $k$?
\end{itemize}
We use Theorem~\ref{AdamsMahdianMahmoodian} to prove that this
problem is also \NP-complete even for bipartite graphs with
maximum degree $4$.
\begin{theorem} {\sc Smallest forcing number of graph} is
\NP-complete for bipartite graphs with maximum degree $4$.
\end{theorem}
\begin{proof}
It is clear that the problem is in \NP. We prove the
\NP-completeness by reducing {\sc Smallest forcing set} to this
probem. Let $G$ be a bipartite graph with maximum degree $3$ and
$M_G$ be a matching in $G$. We construct a new graph $H$ with
maximum degree $4$ as follows:
\begin{itemize}
\item $G$ is a subgraph of $H$, and
\item For any edge $e=\{x,y\} \in E(G)-M_G$, we add vertices $x_e$
 and $y_e$ to $H$ plus three edges $\{x,y_e\}$, $\{x_e,y_e\}$, and $\{x_e,y\}$.
\end{itemize}
Note that $H$ satisfies the conditions of the theorem and any
forcing set for the matching $M_G$ also forces a matching in $H$.
We claim that the smallest forcing number of $H$ is equal to the
smallest forcing number of $M_G$. We can assume that $x_e$ is
matched to $y_e$, otherwise we have the following case: $x_e$ is
matched to $y$ and $y_e$ is matched to $x$. Any forcing set
contains one of these two edges, and choosing one will force the
choice of the other edge. So it is obvious that in this case a
matching 2--switch on these edges will not change the forcing
number. With this assumption, every matching in $G$ corresponds
uniquely to a matching in $H$ and vice versa. For every matching
$M'_G$ in $G$, we denote the corresponding matching in $H$ by
$M'_H$. Now consider a matching $L_G$ in $G$. For every edge $e =
\{x,y\}$ in $L_G - M_G$, the four vertices $x, y, x_e$, and $y_e$
constitute an alternating cycle for $L_H$, so at least one edge
from this alternating cycle should be in the forcing set, and
since choosing the edge $e$ forces the choice of the other edge,
we can assume that $e$ is in the forcing set. Thus a forcing set
$F$ for $L_H$ consists of $L_G - M_G$ plus some edges in $L_G
\cap M_G$. It is not hard to see that $F' = (M_G - L_G) \cup (F
\cap M_G)$ is a forcing set for $M_G$. Since $|M_G - L_G| = |L_G -
M_G|$ and $F \subseteq L_G$, we have $|F'| =|F|$.
\end{proof}

For the problem of finding a smallest forcing set for a given
matching in a planar graph, we have a polynomial
algorithm~\cite{MR99c:05156}, so it is interesting to ask the
following question:
\begin{question}
What is the computational complexity of the following problem:
Given a planar graph $G$, find the smallest forcing number of $G$.
\end{question}
After studying the computational complexity of the problem of
finding the smallest forcing number of a graph it is natural to
do the same for the largest forcing number. So we ask also the
following question, and leave it as an open problem.
\begin{question}
What is the computational complexity of the following problem:
Given a graph $G$, find the largest forcing number of $G$.
\end{question}
\section*{Acknowledgements}
We thank Mohammad Mahdian for many valuable discussions towards
results of this paper, and David de Wit for grammar-bashing.
%
\bibliographystyle{siam}
\bibliography{sr35}
\end{document}